\documentclass[12pt]{article}

    \usepackage{amsmath,amssymb,amsthm}
    \usepackage{graphicx}
    \usepackage{hyperref}
    \usepackage{amsrefs} 

\title{A Discrete Fourier Kernel and Fraenkel's Tiling Conjecture}
\author{
    Ron Graham
        \thanks{
            Department of Computer Science \& Engineering,
            University of California, San Diego.
            Research supported in part by Grant CCR-0310991.}
        \and
    Kevin O'Bryant
        \thanks{
            Department of Mathematics,
            University of California at San Diego, \href{mailto:kevin@member.ams.org}{kevin@member.ams.org}.
            National Science Foundation Mathematical Sciences Postdoctoral Fellow,
            NSF grant DMS-0202460.}
        }
\date{\today}

    \newcommand{\MathReview}[1]{~\href{http://www.ams.org/mathscinet-getitem?mr=#1}{\mbox{\bf MR~#1}}}

    \newcommand{\ZZ}{{\mathbb Z}}
    \newcommand{\QQ}{{\mathbb Q}}
    
    \newcommand{\CC}{{\mathbb C}}
    
    \newcommand{\B}{{\mathcal B}}

    \newcommand{\floor}[1]{\left\lfloor #1 \right\rfloor}
    \newcommand{\ceiling}[1]{\mbox{$\left\lceil #1 \right\rceil$}}
    \newcommand{\fp}[1]{\mbox{$\left\{ #1 \right\}$}}

    \newcommand{\tf}[1]{\mbox{$\left[\!\left[ #1 \right]\!\right]$}}

    \newtheorem{thm}{Theorem}[section]
    \newtheorem{lem}[thm]{Lemma}
    
    \newtheorem{cor}[thm]{Corollary}
    
    \newtheorem{cnj}[thm]{Conjecture}

\begin{document}
\maketitle

\begin{abstract}
The set $\B_{p,r}^q:=\{\floor{nq/p+r} \colon n\in \ZZ\}$ (with integers $p,q,r$) is a Beatty set with density
$p/q$. We derive a formula for the Fourier transform $$\widehat{\B_{p,r}^q}(j):=\sum_{n=1}^p e^{-2\pi ij
\lfloor{nq/p+r}\rfloor/q}.$$ A. S. Fraenkel conjectured that there is essentially one way to partition the
integers into $m\ge 3$ Beatty sets with distinct densities. We conjecture a generalization of this, and use
Fourier methods to prove several special cases of our generalized conjecture.
\end{abstract}



\section{Introduction}

The result that drives this paper is
    \begin{equation}\label{equ:driver}
    \sum_{n=1}^p e^{2\pi i \lfloor{nq/p}\rfloor/q} = \frac{1-e^{2\pi i/q}}{1-e^{2\pi i \bar{p}/q}},
    \end{equation}
where $p$ and $q$ are relatively prime positive integers and $\bar p$ is the multiplicative inverse of
$p\bmod{q}$. The LHS is complicated by the irregularity arising from the floor function, while the RHS is
complicated by the presence of a modular inverse; therein lies the beauty and utility of~\eqref{equ:driver}.

Before stating the general result of which~\eqref{equ:driver} is a special case (Theorem~\ref{thm:MainTheorem}),
we need to introduce some notation. We set $\omega:=e^{2\pi i /q}$, and whenever the range of a summation is not
written explicitly, it is to be taken over all of $\ZZ_q$, the integers modulo $q$:
    $$\sum_{x} = \sum_{x\in \ZZ_q} = \sum_{x=1}^q = \sum_{x=0}^{q-1} .$$
We use the Fourier transform
    $$\hat{f}(j) := \sum_x f(x) \omega^{-jx}$$
($\hat f(j)$ is called the $j$-th Fourier coefficient), the Fourier inversion formula
    $$\hat{\hat{f}}(x) = q f(x),$$
convolution
    $$f\ast g (x) := \sum_y f(y) {g(x-y)}$$
and the interchange-of-summations result
    $$\widehat{f\ast g}(j) = \hat{f}(j) \hat{g}(j).$$

Also, let
    $$\tf{P}\,R := \begin{cases} 0 & \text{$P$ is False;}\\ R & \text{$P$ is True.}\end{cases}$$
Note that $\tf{\text{False}}R$ is defined even if $R$ is not. When $R=1$, we omit it from the notation. We also
adopt the common practice of identifying a multiset with its indicator function, i.e., $S(x)$ is the
multiplicity of $x$ in the multiset $S$.

We distinguish the rational Beatty sets ($p, q$ are any integers, and $r$ any real number)
    $$\B_{p,r}^q := \left\{ \floor{n \frac qp+r}  \colon n \in \ZZ \right\}.$$
Usually, $q$ will be fixed and in this situation we omit it from the notation. We will always assume that $r$ is
an integer\footnote{There is no loss of generality in assuming that $r$ is an integer. To see this, let $n_0$ be
a value of $n$ for which the fractional part $\fp{n q/p+r}$ is minimal; then $\B_{p,r}^q = \B_{p,\floor{n_0
q/p+r}}^q$.}, and when $r=0$ we omit it from our notation. Note that the density of the set $\B_{p,r}^q$ is
$p/q$.

Note that $\B_{p,r}^q$ consists of $p$ congruence classes modulo $q$, and so $\B_{p,r}^q$ is naturally
considered as a subset of $\ZZ_q$. If the $p$ points were perfectly evenly distributed around $\ZZ_q$ (as
happens if $q/p\in\ZZ$), then the Fourier transform would be 0 except at multiples of the difference between
points. Thus, one naturally expects that $\B_{p}(j)$ will be small except when $jq/p$ is near an integer. This
is confirmed by Figure~1, which shows the Fourier transform of $\B_{24}^{121}$, and Figures~2 and~3, which show
$|\widehat{\B_{p,r}^q}(1)|$ for small relatively prime $p,q$ ($r$ is irrelevant). Theorem~\ref{thm:MainTheorem}
gives an explicit formula for $\widehat{\B_{p,r}^q}$ which quantifies the validity of this expectation.

\begin{figure}
\begin{center}
    \begin{picture}(216,450)
        \put(0,225){\includegraphics{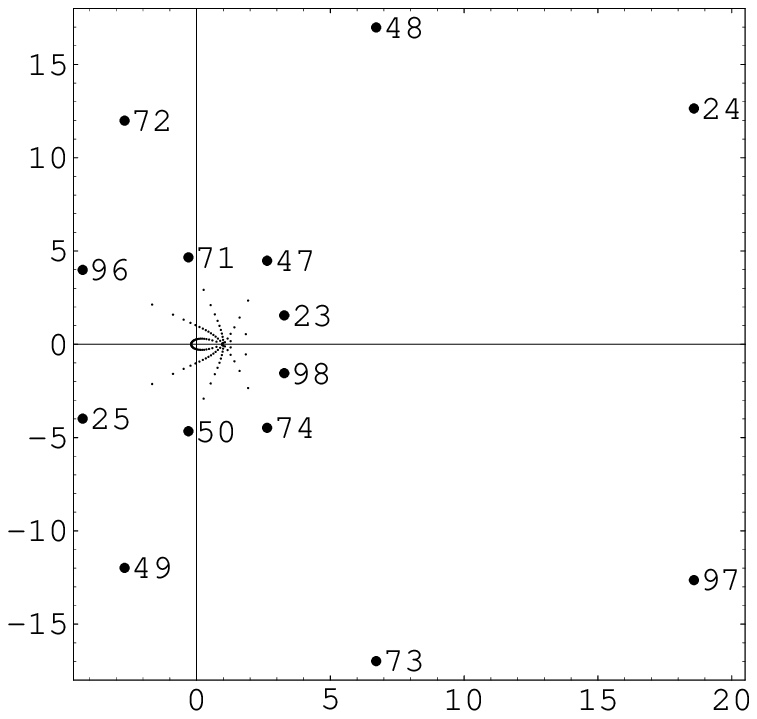}}
        \put(0,0){\includegraphics{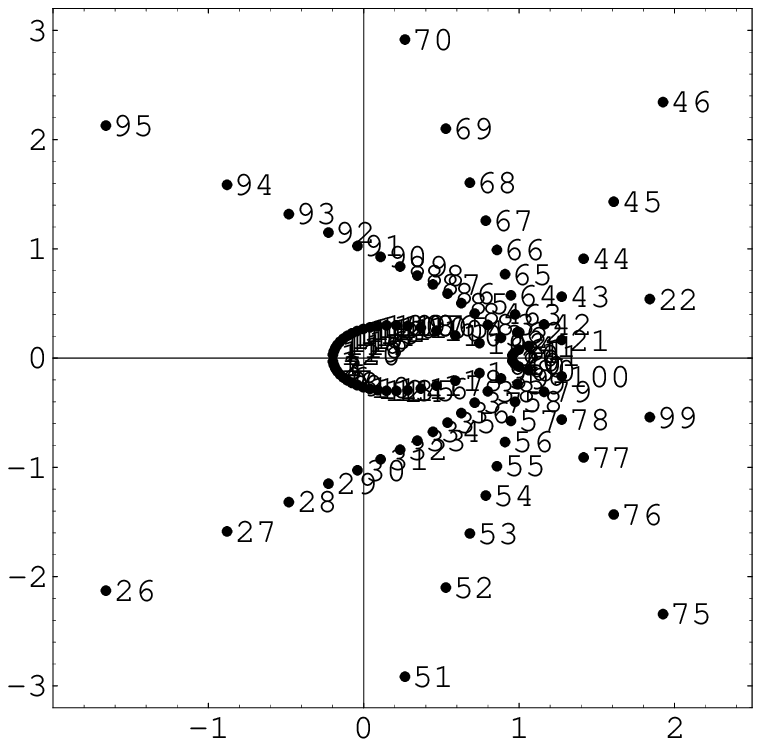}}
    \end{picture}
\end{center}
\caption{The points $\widehat{\B_{24}^{121}}(j)$ ($1\le j \le 120$), shown in the complex plane and labeled by
$j$. The top graph shows all 120 points. The bottom graph, which resembles a spider, shows only those closest to
1.}
\end{figure}

\begin{thm}\label{thm:MainTheorem}
Let $p\not=0$, $q>1$ be integers with $g:=\gcd(p,q)$, let $\bar{p}$ satisfy $p\bar{p} \equiv g \pmod{q}$, and
let $r$ be any integer. Then
    $$
    \widehat{\B_{p,r}^q}(0)=p
    $$
and for $j\not\equiv 0 \pmod{q}$
    \begin{align*}
    \widehat{\B_{p,r}^q}(j)
        &=  \tf{g|j}\, g\, \frac{1-\omega^j}{1-\omega^{j\bar p}} \, \omega^{-jr},\\
    \left|\widehat{\B_{p,r}^q}(j) \right|
        &=  \tf{g|j}\, g\, \left|\frac{\sin (\pi j /q)}{\sin (\pi j \bar p /q)}\right|. \notag
    \end{align*}
\end{thm}

\begin{figure}
\begin{picture}(360,120)
    \put(0,0){\includegraphics{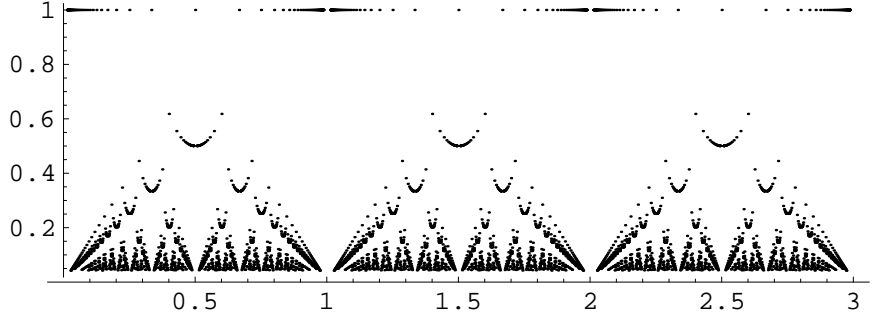}}
\end{picture}
\caption{The points $(\frac pq, |\widehat{\B_p^q}(1)|)$ for $\gcd(p,q)=1$, $0\!<\!p\!<\!3q$, and $1\le q \le
75$}
\end{figure}

\begin{figure}\label{fig:firstcoeff}
\begin{picture}(360,120)
 \put(0,0){\includegraphics{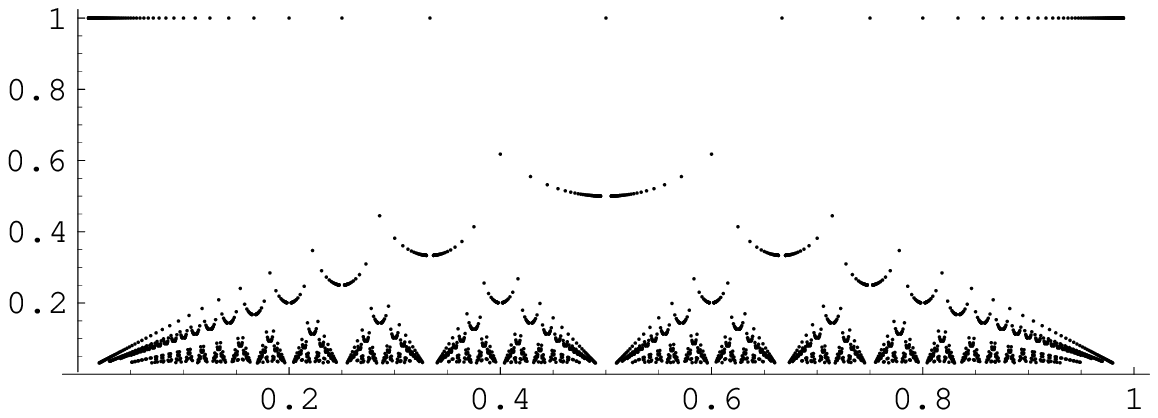}}
\end{picture}
\caption{The points $(\frac pq, |\widehat{\B_{p}^q}(1)|)$, for $\gcd(p,q)=1$, $0< p <q$, $1<q\le 100$.}
\end{figure}

\begin{figure}\label{fig:secondcoeff}
\begin{picture}(360,120)
 \put(0,0){\includegraphics{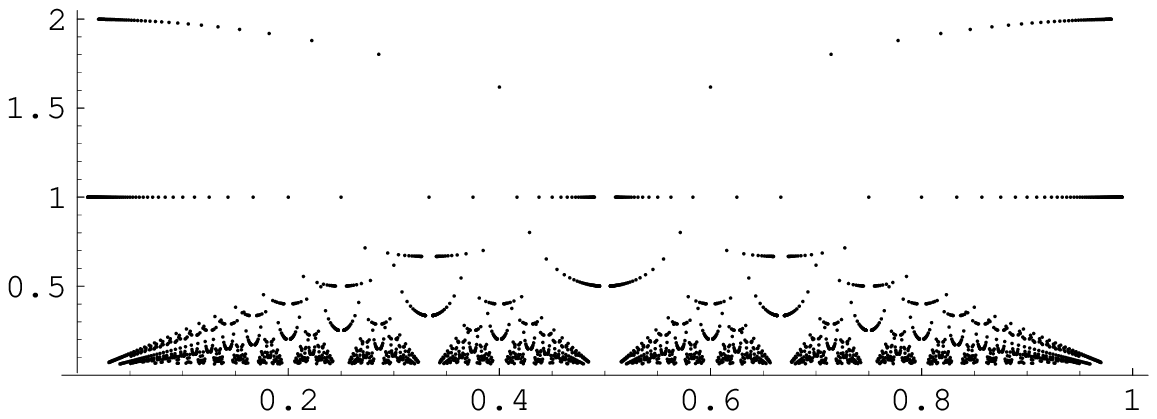}}
\end{picture}
\caption{The points $(\frac pq, |\widehat{\B_{p}^q}(2)|)$, for $\gcd(p,q)=1$, $0< p <q$, $2<q\le 100$.}
\end{figure}

Figure~2 shows $|\widehat{\B_p^q}(1)|$ for relatively prime $p$ and $q$ with $0<p<3q$ and $1<q\le 75$. Figures~3
and~4 show the first and second coefficients when $p$ and $q$ are relatively prime, $0<p<q$, and $q\le 100$.
There are three symmetries visible to the naked eye. The first (from Figure~2) is that $|\widehat{\B_{p+q}^q}| =
|\widehat{\B_{p}^q}|$. In other words, the function $|\widehat{\B_x^q}(j)|$ is periodic with period $q$. This is
a consequence of the fact $\B_{p+q}^q(x)=1+\B_p^q(x)$, which we prove along the way to proving
Theorem~\ref{thm:MainTheorem}.

The second is that
    $$\left| \widehat{\B_{p}^q}(j) \right| = \left| \widehat{\B_{q-p}^q}(j) \right|,$$
which is seen in the pictures as a symmetry about $1/2$. This is a consequence of a theorem of Fraenkel
(Corollary~\ref{cor:FPTrational} below) which states that the complement of a rational Beatty set is a rational
Beatty set. We give a new proof of this in Section~\ref{sec:FPT}.

The third symmetry is that the graphs on $[\frac 14, \frac 13]$, on $[\frac 13, \frac 12]$, etc., seem to be
quite similar. This is essentially the symmetry of the continued-fraction map $x\mapsto 1/x\pmod{1}$. For each
rational $p/q$, there is a unique finite sequence $[a_0;a_1,\dots,a_n]$ of integers with the properties: for
$i>0$, $|a_i|\ge 2$; for $0<i<n$, if $a_i=\pm2$ then $a_ia_{i+1}$ is positive; $a_n\not= -2$; and
    $$\frac pq = a_0+\frac{1}{a_1+\frac{1}{a_2+ \ddots +\frac{1}{a_n}}}.$$
This is the nearest-integer continued fraction (commonly abbreviated NICF). Compare Figure~3 with Figure~5,
which shows the points $(\frac pq, \prod_{i=1}^n |a_i| )$: the points in Figure~5 are located precisely at the
bottom of the ``cups'' in Figure~3.

We remark that while Eq.~\eqref{equ:driver} connects Beatty sequences with density $p/q$ directly to the inverse
of $p$ modulo $q$, the direct connection between both objects and continued fractions is well-studied. We also
note that the inverse of $p$ modulo $q$ has arisen independently in the recent work of
Simpson~\cite{2004.Simpson}, in which he uses generating functions to prove necessary and sufficient conditions
for two rational Beatty sequences to not intersect.

\begin{figure}
\begin{picture}(360,120)
    \put(0,0){\includegraphics{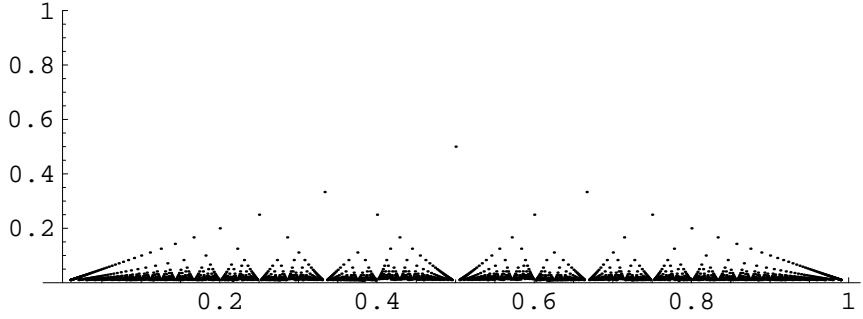}}
\end{picture}
\caption{The points $(p/q,\prod_{i=1}^n |a_i|)$, where $[a_0;a_1,\dots,a_n]$ is the NICF of $p/q$}
\end{figure}

The main result in the study of Beatty sequences was discovered by Lord Rayleigh: If $\alpha$ is irrational and
$\frac 1\alpha+\frac1\beta=1$, then the sets $\{\floor{n\alpha}\colon n\ge 1\}$, $\{\floor{n\beta}\colon n\ge
1\}$ partition $\ZZ^+$. In the 1950s Skolem extended to this non-homogeneous sets, and in 1969 Fraenkel
corrected Skolem's work and extended it to include rational $\alpha$. We direct the reader
to~\cite{MR2004g:11017} for the general rational/irrational statement, an elementary proof, and the history of
Fraenkel's Partition Theorem. The rational case (Corollary~\ref{cor:FPTrational}) is an easy consequence of
Theorem~\ref{thm:MainTheorem}.

\begin{cor}[Fraenkel's Partition Theorem, rational case]\label{cor:FPTrational}
The sets $\B_{p_1,r_1}^q$ and $\B_{p_2,r_2}^q$, where $q,p_k,r_k$ are integers, partition $\ZZ$ if and only if
$p_1+p_2=q$ and $$p_1r_1+p_2r_2\equiv -\gcd(p_1,q) \pmod{q}.$$
\end{cor}

Attempts to extend this to more than two sequences have had some success, but a general statement remains
elusive. In the early 1970s (see \cites{MR82j:10001,MR46:8875}), Fraenkel was led to conjecture that there was
essentially only one way to partition $\ZZ$ into Beatty sets with distinct densities.

\begin{cnj}[Fraenkel's Conjecture]
If the sets
    $$\left\{ \floor{n\alpha_1+r_1} \colon n \in \ZZ \right\},
        \left\{ \floor{n\alpha_2+r_2} \colon n \in \ZZ \right\}, \ldots,
        \left\{ \floor{n\alpha_m+r_m} \colon n \in \ZZ \right\}$$
partition $\ZZ$, $m\ge 3$, and $1<\alpha_1<\alpha_2< \dots < \alpha_m$, then $\alpha_k=2^k-2^{k-m}$.
\end{cnj}

The conjecture has been proven in case any $\alpha$ is irrational~\cite{MR48:3911}, $\alpha_1\le
1.5$~\cite{MR93b:11023}, or $m\le 6$~\cite{MR2001f:11039}, and in several other less-easily-stated
circumstances. The article of R. Tijdeman~\cite{MR2001h:11011} contains an excellent survey of the progress on
Fraenkel's Conjecture.

We offer a stronger conjecture, and will apply Theorem~\ref{thm:MainTheorem} to prove some special cases. We say
that sets $S_1, \dots, S_m$ are a {\em perfect $c$-fold covering of $\ZZ$} if $S_1(x)+S_2(x)+\dots+S_m(x)=c$ for
all $x\in \ZZ$, and simply that they are a {\em perfect covering} if we do not wish to specify $c$.

\begin{cnj}[Covering Fraenkel Conjecture]
Let $q,m\ge 3$, and let $p_1,p_2,\dots,p_m$ be distinct integers with $0<p_k<q$,
\mbox{$\gcd(q,p_1,\dots,p_m)=1$}, and with no proper subset $I\subsetneq[m]$ having $\sum_{i\in I} p_i\equiv 0
\pmod{q}$. Let $r_1,\dots,r_m$ be arbitrary integers. The sets $\B_{p_k,r_k}^q$ ($1\le k \le m$) are a perfect
covering if and only if:
    \renewcommand{\theenumi}{\roman{enumi}}
    \begin{enumerate}
    \item[(i)]   $q$ is odd, and $m$ is the order of 2 modulo $q$;
    \item[(ii)]   the sets can be renumbered so that there are integers $\delta$ (relatively prime to $q$ with modular
            inverse $\bar \delta$)
            and $\gamma$ such that
            for $1\le k \le m$:
            $$p_k \equiv \delta 2^{m-k} \pmod{q}
                \quad \text{and}\quad
                r_k \equiv \gamma - \bar \delta 2^{k-1}\pmod{q}.$$
    \end{enumerate}
\end{cnj}

This is admittedly not as pithy as Fraenkel's Conjecture, but we hope that its greater generality will shed new
light on an old problem. Since $\B_{p,r}^q(x)=1+\B_{p+q,r}^q(x)$, we may assume without loss of generality that
$p_k<q$. We use Theorem~\ref{thm:MainTheorem} to prove the ``if'' claim of the CFC, and in fact we show that if
the covering (with all $p_k<q$) is a perfect $c$-fold covering, then $c$ is the number of ones in the binary
expansion of $\delta(2^m-1)/q$. We prove the ``only if'' part of the CFC under the additional hypothesis that
$m\le 5$, and also under the additional hypotheses that $\gcd(q,p_k)=1$ for some $k$ and that $q$ is
sufficiently large.

We suspect that if $m\ge 3$ sequences $S(\alpha_k,\beta_k)=\{\floor{n\alpha_k+\beta_k}\colon n\in\ZZ\}$ are a
perfect covering and some $\alpha_k$ is irrational, then two of the $\alpha_k$'s are in fact equal. Graham's
proof \cite{MR48:3911} of this in the 1-covering case does not extend easily to multiple coverings. We have not
investigated this further, and at any rate this guess does not fall within the scope of the present article.

We will use the following consequence of Theorem~\ref{thm:MainTheorem} several times.
\begin{cor}[Covering Criterion]\label{cor:CoveringCriterion}
The sequences $\B_{p_k,r_k}^q$ ($1\le k \le m$) are a perfect $c$-fold covering of $\ZZ$ if and only if
    \begin{equation*}
    cq = \sum_{k=1}^m p_k
    \end{equation*}
and for $1\le j < q$ (with $g_k=\gcd(p_k,q)$),
    \begin{equation*}
    0 = \sum_{k=1}^m \tf{g_k | j} g_k \frac{\omega^{-jr_k}}{1-\omega^{j \bar p_k}}.
    \end{equation*}
\end{cor}

\section{Proof of Theorem~\ref{thm:MainTheorem}}

We will give a kernel function $K(x)$, whose transform we can compute, such that
    \begin{equation} \label{equ:convplan}
    \B_{p,r}^q \ast K(x) = R(x),
    \end{equation}
where $R(x)$ will also have an easily computed transform. Taking the Fourier transform of this equation gives
$\widehat{\B_{p,r}^q}(j) \hat{K}(j) = \hat{R}(j)$, which is the same as $\widehat{\B_{p,r}^q}(j)=
\hat{R}(j)(\hat{K}(j))^{-1}$. The zeroth coefficient can be dealt with immediately, and thereafter we
demonstrate that we may assume that $r=0$, that $\gcd(p,q)=1$, and that $0<p<q$. Then, we define $K$ and $R$,
show that~\eqref{equ:convplan} holds, and compute the transforms of $R$ and $K$ to complete the proof.

We have
    $$\widehat{\B_{p,r}^q}(j) := \sum_x \B_{p,r}^q(x) \omega^{-j x} = \sum_{n=0}^{p-1}
    \omega^{-j\floor{nq/p+r}}.$$
For $j=0$, this gives the first claim of Theorem~\ref{thm:MainTheorem}: $\widehat{\B_{p,r}^q}(0)=p$.

Since we assume that $r$ is an integer, we also have
    $$\widehat{\B_{p,r}^q}(j) =\sum_{n=0}^{p-1}
    \omega^{-j\floor{nq/p+r}} = \sum_{n=0}^{p-1} \omega^{-j\floor{nq/p}-jr}
    =\omega^{-jr}\widehat{\B_{p,0}^q}(j).$$
Thus, it is sufficient to work with $r=0$.

Now we wish to show that we may take $p$ and $q$ to be relatively prime. If $S\colon \ZZ\to\CC$ is periodic with
period $\ell$, and $S_\ell$ and $S_{k\ell}$ are the induced functions on $\ZZ_\ell$ and $\ZZ_{k\ell}$, then
    $$\widehat{S_{k\ell}}(j)= k \tf{k|j} \widehat{S_\ell}(j/k).$$
Thus, if $ga=p$ and $gb=q$ with $\gcd(a,b)=1$, then
    $$\widehat{\B_p^q}(j) = g \tf{g|j} \widehat{\B_a^b}(j/g).$$
If we assume for the moment that we have proved Theorem~\ref{thm:MainTheorem} in the relatively prime case, then
we have
    \begin{align*}
    \widehat{\B_p^q}(j) &= g \tf{g|j} \frac{1-(e^{2\pi i /b})^{j/g}}{1-(e^{2\pi i \bar{a} /b})^{j/g}} \\
                    &=g \tf{g|j} \frac{1-(e^{2\pi i /(bg)})^{j}}{1-(e^{2\pi i \bar{a} /(bg)})^{j}}\\
                    &=g \tf{g|j} \frac{1-\omega^j}{1-\omega^{j\bar{a}}}
    \end{align*}
where $\omega=e^{2\pi i /q}$ and $\bar{a}$, the inverse of $a=p/g$ modulo $b=q/g$ satisfies $p\bar{a}\equiv g
\pmod{q}$. Thus, it is sufficient to work with relatively prime $p$ and $q$.

We may assume that $0<p<q$ since
    \begin{equation}\label{equ:density.is.periodic}
    \B_{p+q}^q(x)=1+\B_p^q(x),
    \end{equation}
and so the Fourier coefficients (except the zeroth) are sensitive only to $p\bmod{q}$. To
see~\eqref{equ:density.is.periodic}, observe that the function $\B_p^q(x)$ takes on only the values
$\floor{p/q}$ and $\ceiling{p/q}$. Moreover, $\B_{p}^q(x)=\beta$ (with $\beta = \ceiling{p/q}$) if and only if
there is an integer $n$ with $x\le nq/p < (n+\beta-1)q/p < x+1$, which is equivalent to
    $$x \frac pq \le n <  (x+1) \frac pq - (\beta-1) = x\frac pq + \frac pq-\beta+1.$$
This happens if and only if the fractional part $\fp{xp/q}$ is 0 or is strictly larger than
$\beta-p/q=1-\fp{p/q}$. Thus, the property $\B_p^q(x)=\ceiling{p/q}$ can be described entirely in terms of the
fractional part $\fp{x p/q}$, which depends only on $p\bmod q$.

Now, we take $K$ to be the set $\{1-\bar{p},2-\bar{p},\dots, 0\}=(-\bar p,0]$, and set $R(x):= \left| \B_p \cap
[x,x+\bar p) \right|$.

We use the following two properties of the Beatty set $\B_p$ (with $0<p<q$). We call the first property
``duality'': an integer $k\in \B_p$ if and only if the fractional part of ${k p/q}$ is 0 or strictly greater
than $1-p/q$. To prove this simply observe that $k\in \B_p$ if there is an integer $n$ with $k\le n q/p < k+1$,
which we rearrange as $kp/q \le n < kp/q + p/q$. This happens exactly if the fractional part of $kp/q$ is 0 or
strictly greater than $1-p/q$. The second property is called ``balance'': for all real numbers $x<y$, the number
$\left|\B_p \cap [x,y) \right|$ is either $\floor{(y-x)p/q}$ or $\ceiling{(y-x)p/q}$. To prove this, observe
that we are counting the integers $n$ with $x\le n q/p <y$, which is equivalent to $xp/q \le n <yp/q$. Since we
only care about integral $n$, we can write this as $\ceiling{x p /q }\le n <\ceiling{yp/q}$. Clearly there are
$\ceiling{yp/q}-\ceiling{x p /q }$ such $n$, and this is
    $$\left\lceil \dfrac{yp}{q} \right\rceil - \left\lceil \dfrac{xp}{q} \right\rceil
        =\left(\frac{yp}{q}+\epsilon_1\right)-\left(\frac{xp}{q}+\epsilon_2\right)
        =(y-x)\frac pq +\epsilon_1-\epsilon_2,$$
where the `$\epsilon$'s are both in $[0,1)$. It follows that the true value is an integer which is strictly
less than 1 away from $(y-x)p/q$, i.e., either $\floor{(y-x)p/q}$ or $\ceiling{(y-x)p/q}$.

We will show that
    \begin{equation}\label{equ:R}
        R(x) = \floor{\bar p p/q} + \tf{x=0},
    \end{equation}
but first we show that $R(0)>R(1)$. Since $R(0)$ counts the number of elements of $\B_p=\{\floor{nq/p} \bmod q
\colon 1\le n \le p\}$ in $[0,\bar p)$ and $R(1)$ counts those in $[1,\bar p]$, we need to show that $0\in \B_p$
and $\bar p \not\in \B_p$. Obviously $\floor{p q/p}\equiv  0 \bmod q$, so the substance here is that $\bar p
\not\in \B_p$. By duality, $\bar p \in \B_p$ if and only if the fractional part of $\bar p p/q$ is 0 or strictly
greater than $1-p/q$. Since $p\bar p \equiv 1\pmod q$, the fractional part of $\bar p p /q $ is $1/q$, which is
neither 0 nor strictly greater than $1-p/q$ (using that $p,q$ are relatively prime).

We now show~\eqref{equ:R} by evaluating $\sum_x R(x)$ in two ways. First, every $y\in\B_p$ contributes to $R(y),
R(y-1), \dots, R(y-\bar p+1)$. Thus
    \begin{equation}\label{equ:sumR1}
    \sum_x R(x) = |\B_p| \bar p = p \bar p.
    \end{equation}
Second, by the balance property of Beatty sets, we know that $R(x)$ is either $\floor{\bar p p/q}$ or
$\ceiling{\bar p p/q}=\floor{\bar p p/q}+1$, and in particular $\left|R^{-1}(\floor{\bar p
p/q})\right|+\left|R^{-1}(\floor{\bar p p/q}+1)\right|=q$. Thus
    \begin{align}
    \sum_x R(x) &= \left|R^{-1}(\floor{\bar p p/q})\right| \floor{\bar p p/q} +
                        \left|R^{-1}(\floor{\bar p p/q}+1)\right| \left(\floor{\bar p p/q}+1\right) \notag\\
                &=  q \floor{\bar p p/q}+\left|R^{-1}(\floor{\bar p p/q}+1)\right|. \label{equ:sumR2}
    \end{align}
Reducing~\eqref{equ:sumR1} and~\eqref{equ:sumR2} modulo $q$ tells us that $\left|R^{-1}(\floor{\bar p
p/q}+1)\right|=1$. Since $R(0)>R(1)$, we know that $R^{-1}(\floor{\bar p p/q}+1)=\{0\}$, whence
Eq.~\eqref{equ:R}.

Theorem~\ref{thm:MainTheorem} now follows from the straightforward calculations (for $j\not\equiv0\pmod{q}$)
    \begin{align*}
    \hat{K}(j) &= \frac{1-\omega^{j\bar p}}{1-\omega^j} \\
    \hat{R}(j) &= 1
    \end{align*}
and
    \begin{align*}
    \B_{p}\ast K(x) &:=  \sum_{y} \B_{p}(y) {K(x-y)} \\
    &= \, \sum_{y}\tf{y \in \B_{p}} \tf{-\bar{p} < x-y \le 0} \\
    &= \, \sum_{y}\tf{y \in \B_{p}} \tf{y\in [x,x+\bar{p})} \\
    &= \, \left|\B_{p} \cap [x,x+\bar{p})\right|
    &=: R(x).
    \end{align*}

\subsection{An Interesting Variation}

There is another interesting way to finish the proof of Theorem~\ref{thm:MainTheorem}. It plays on another
expression of the ``balance'' property of Beatty sets: for fixed $t$, the difference
$\floor{(n+t)q/p}-\floor{nq/p}$ is either $\floor{tq/p}$ or $\ceiling{tq/p}$.

Assume that $j\not\equiv0\pmod{q}$, $r=0$, and $p<q$ are relatively prime and positive.

Let $\bar q$ be the inverse of $q$ modulo $p$, and let $\bar p$ be the inverse of $p$ modulo $q$. We will use
the identity
    $$\floor{q\bar q/p} = (q\bar q -1)/p \equiv -\bar p \pmod{q}.$$
Set $b(n):=\floor{nq/p}$, and
    $$\Delta(n) :=b(n)-b(n-\bar q).$$
We compute $\sum_{n=0}^{p-1} \Delta(n)$ in two ways. First, the sum telescopes to
    \begin{multline*}
    \sum_{n=0}^{p-1}\Delta(n) = \sum_{n=p-\bar q}^{p-1} b(n) - \sum_{n=-\bar q}^{-1} b(n)
            = \sum_{n=-\bar q}^{-1} b(n+p)- \sum_{n=-\bar q}^{-1} b(n)\\
            =\sum_{n=-\bar q}^{-1} \big(q+b(n)\big)- \sum_{n=-\bar q}^{-1} b(n) =
            \bar q q.
    \end{multline*}
Second, note that $\Delta(n)$ is either $a:=\floor{\bar qq/p}$ or $a+1=\ceiling{\bar q q/p}$, say there are
$\beta$ integers $n$ inclusively between 0 and $p-1$ with $\Delta(n)=a+1$, and $p-\beta$ integers $n$ with
$\Delta(n)=a$. We have
    $$\sum_{n=0}^{p-1} \Delta(n) = \beta(a+1)+(p-\beta)a=\beta+pa.$$
Equating these two evaluations modulo $p$ (and using $0\le \beta\le p$), we find that $\beta=1$. By direct
arithmetic, $\Delta(0)=\ceiling{\bar qq/p}$, and so for $x\not\equiv0\pmod{p}$, $\Delta(x)=\floor{\bar qq/p}$.

We now use this information directly (set $\gamma=\omega^{-j}$):
    \begin{align*}
    \hat{\B_p}(j)   &=  \sum_{n=1}^p \gamma^{\floor{nq/p}} = 1+\sum_{n=1}^{p-1} \gamma^{b(n)}\\
        &=  1+ \sum_{n=1}^{p-1} \gamma^{b(n)}\gamma^{b(n-\bar q)-b(n)}\gamma^{\floor{\bar q q/p}} \\
        &=  1+ \gamma^{\floor{\bar q q/p}} \sum_{n=1}^{p-1} \gamma^{b(n-\bar q)}\\
        &=  1+ \gamma^{\floor{\bar q q/p}}
                    \left(\Big(\sum_{n=1}^{p} \gamma^{b(n-\bar q)}\Big) - \gamma^{b(p-\bar q)}\right)\\
        &=  1+ \gamma^{\floor{\bar q q/p}}\left( \hat{\B_p}(j) - \gamma^{-\floor{\bar q q/p}-1}\right)
    \end{align*}
Solving this equation yields
    $$\hat{\B_p}(j) = \frac{1-\gamma^{-1}}{1-\gamma^{\floor{\bar q q /p}}} =
    \frac{1-\omega^j}{1-\omega^{-j\floor{\bar q q/p}}}= \frac{1-\omega^j}{1-\omega^{j\bar p}}.$$

\section{Proof of Fraenkel's Partition Theorem}\label{sec:FPT}

By Corollary~\ref{cor:CoveringCriterion}, we can assume that $p_1+p_2=q$; we need only show that
    \begin{equation}
    0 = \tf{g_1 | j} \,g_1\, \frac{\omega^{-jr_1}}{1-\omega^{j \bar p_1}}\;
            +\;\tf{g_2 | j} \,g_2 \,\frac{\omega^{-jr_2}}{1-\omega^{j \bar p_2 }}\label{equ:2pieces}
    \end{equation}
is satisfied for $1\le j<q$ if and only if $p_1r_1+p_2r_2\equiv -g_1 \pmod q$.

We first assume that~\eqref{equ:2pieces} holds for all $j\in[1,q)$. In particular, we set $j=p_1$. Since
$p_1+p_2=q$, we have $g_1=g_2,\, g_1|j$, and~\eqref{equ:2pieces} simplifies to
    \begin{equation*}
    \frac{\omega^{-p_1r_1}}{1-\omega^{p_1\bar p_1}}  +
        \frac{\omega^{-p_1r_2}}{1-\omega^{p_1\bar p_2}} =0.
    \end{equation*}
Rearranging this gives
    \begin{equation}\label{equ:twoterms}
    \frac{\omega^{-p_1r_2}}{\omega^{-p_1r_1}}
        = - \frac{1-\omega^{p_1 \bar p_2}}{1-\omega^{p_1\bar p_1 }}.
    \end{equation}
Since $-p_1\equiv p_2 \pmod q$, $p_1 \bar p_1 \equiv g_1 \pmod q$, and $p_1 \bar p_2  \equiv -p_1 \bar p_1 \pmod
q$, Eq.~\eqref{equ:twoterms} becomes
    $${\omega^{p_2r_2+p_1r_1}}
        = \frac{\omega^{-p_1r_2}}{\omega^{-p_1r_1}}
        = -\frac{1-\omega^{-g_1}}{1-\omega^{g_1}}
        = \omega^{-g_1},$$
whence $p_2r_2+p_1r_1\equiv -g_1 \pmod q$. We can read this argument from the bottom up to see the other half of
``if and only if''.

\section{Fraenkel's Covering Conjecture}

\subsection{Constructions}
\begin{lem}\label{lem:constructions}
If $q\ge 3$, $2^m \equiv 1 \pmod{q}$, $p_k\equiv  2^{m-k} \pmod{q}$, and $r_k\equiv -2^{k-1}\pmod{q}$, then
$\B_{p_k,r_k}$ ($1\le k \le m$) is a perfect covering.
\end{lem}

\begin{proof}
We have
    $$\sum_{k=1}^m p_k \equiv \sum_{k=1}^m 2^{m-k} =2^m-1 \equiv 0 \pmod{q}$$
so the first equation of the Covering Criterion is satisfied for some $c$. Our hypotheses imply that
$\gcd(p_k,q)=1$ and $\bar p_k\equiv 2^k \pmod{q}$, so we need only to show that
    \begin{equation}
    0=\sum_{k=1}^m \frac{\omega^{-r_k}}{1-\omega^{\bar p_k}}
    = -\sum_{k=1}^m \frac{\omega^{2^{k-1}}}{1-\omega^{2^k}} \label{equ:tempsum}
    \end{equation}
holds for $\omega$ any $q$-th root of unity except 1.

Since $1-\omega^{2^k} = (1-\omega)\prod_{s=0}^{k-1} (1+\omega^{2^s})$, we can bring the summands
in~\eqref{equ:tempsum} over a common denominator:
    $$\sum_{k=1}^m \frac{\omega^{2^{k-1}}}{1-\omega^{2^k}}
        =
        \sum_{k=1}^m \frac{\omega^{2^{k-1}} \prod_{s=k}^{m-1} (1+\omega^{2^s})}
                          {(1-\omega) \prod_{s=0}^{m-1} (1+\omega^{2^s})},$$
and we see that it will suffice to show that
    $$\sum_{k=1}^m {\omega^{2^{k-1}} \prod_{s=k}^{m-1} (1+\omega^{2^s})}$$
is zero. But
    $${\omega^{2^{k-1}} \prod_{s=k}^{m-1} (1+\omega^{2^s})}=\sum_{a\in A_k} \omega^a,$$
where $A_k$ consists of those integers whose binary expansions have the form $(b_{m-1}b_{m-2} \cdots b_1 b_0)_2$
with $b_0=b_1 = \dots b_{k-2}=0$ and $b_{k-1}=1$. Thus
    $$\sum_{k=1}^m {\omega^{2^{k-1}} \prod_{s=k}^{m-1} (1+\omega^{2^s})} =
        \sum_{k=1}^m \sum_{a\in A_k} \omega^a = \sum_{x=1}^{2^m-1} \omega^x.$$
and since $2^m-1$ is a multiple of $q$, this is zero.
\end{proof}

\begin{thm}\label{thm:Multiplicity}
Suppose that $q,m\ge3,p_k,r_k,\delta,\gamma$ ($1\le k \le m$) satisfy conditions (i) and (ii) of the CFC, and
that that $0<p_k<q$ for $1\le k \le m$. The $m$ sequences $\B_{p_k,r_k}^q$ are a perfect $c$-fold covering,
where $c$ is the number of ones in the binary expansion of $\delta (2^m-1)/q$.
\end{thm}

\begin{proof}
This is equivalent to something by the Covering Criterion: the zeroth coefficient criterion becomes
$\sum_{k=1}^m p_k = cq$, and the non-zero coefficient part of the criterion, with an appropriate choice of $j$
and multiplying by $\omega^{-j\gamma}$, becomes equivalent to Lemma~\ref{lem:constructions}. Thus, what remains
to be proved is that $\sum_{k=1}^m p_k =cq$.

We have
    $$S:=\sum_{k=1}^m p_k \equiv \sum_{k=1}^m \delta 2^{m-k} = \delta (2^m-1) \equiv 0 \pmod q$$
so $S$ is definitely a multiple of $q$. Note that $\zeta_k$ defined by
    \[
    \zeta_k\equiv 2^{m-k}\delta (2^m-1)/q \pmod{2^m-1}
    \qquad \text{and} \qquad
    0<\zeta_k<2^m
    \]
also satisfies $\zeta_k=p_k(2^m-1)/q$. If the binary expansion of $\zeta_1$ is $(b_{m-1}b_{m-2}\cdots
b_1b_0)_2$, then the binary expansion of $\zeta_k$ is $(b_{k-2}b_{k-3}\cdots b_0 b_{m-1} \cdots b_{k-1})_2$. It
follows that $\sum_{k=1}^m \zeta_k = (2^m-1) w$, where $w$ is the number of `1's in the binary expansion of any
of the $\zeta_k$, in particular, in the expansion of $\zeta_m=\delta (2^m-1)/q$. On the other hand,
$\sum_{k=1}^m \zeta_k = \sum_{k=1}^m p_k (2^m-1)/q= c (2^m-1)$. Consequently, $c=w$.
\end{proof}

According to Fraenkel's Conjecture, the values of $q$ for which there is a nontrivial perfect 1-covering by $\ge
3$ Beatty sets are of the form $2^m-1$. As a consequence of the preceding result, we can identify those $q$
which allow for a nontrivial perfect 2-covering.

\begin{cor}
If $q$ is of the form $2^m-1, m \ge 3$, or $\frac{2^{2uv}-1}{2^u+1}$ for $u,v \ge 1$, then there is a perfect
$2$-covering by at least three Beatty sets with period $q$ and distinct densities.
\end{cor}

\begin{proof}
In this proof, we abbreviate $\gcd(a,b)$ as simply $(a,b)$. From Theorem~\ref{thm:Multiplicity}, we know that
$q$ has a perfect $2$-covering if there is a $\delta < q$ with $(\delta,q)=1$ such that $\delta(2^m-1)/q =
2^s+1$ for some $s < m$, where $m$ is the order of $2$ modulo $q$. Let $d = (2^m-1,2^s+1)$. Since the fractions
    $$\frac{(2^m-1)/d}{(2^s+1)/d} = \frac{q}{\delta}$$
are both reduced we see that
    $$q=\frac{2^m-1}{(2^m-1,2^s+1)}.$$
Using the elementary identity $2^{(a,b)}-1=(2^a-1,2^b-1)$, we now have
    \begin{align*}
    2^{(m,2s)}-1&=(2^m-1,2^{2s}-1)\\
    &=(2^m-1,2^s+1)\,(2^m-1,2^s-1)\\
    &= (2^m-1,2^s+1)\; (2^{(m,s)}-1).
    \end{align*}
There are two possibilities: $(m,2s)=(m,s)$ or $(m,2s)=2\,(m,s)$. In the first case, we find that
$(2^m-1,2^s+1)=1$ and so $q=2^m-1$. It is clear that any value of $m \ge 3$ will work here. In the second case,
we must have $m=2M$. Thus, in this case $(m,2s)=(2M,2s)=2\,(M,s)=2\,(m,s)$ by hypothesis, i.e., $(M,s)=(2M,s)$.
Hence, if we let $u=(M,s)$ then $M=uv$ for some $v$. This implies that
    $$q=\frac{2^{2uv}-1}{2^u+1}$$
for some $u,v \ge 1$. This completes the proof.
\end{proof}

If our Covering Fraenkel Conjecture is correct, then these are all the values of $q$ for which there exist
non-trivial perfect $2$-coverings. Besides the values of the form $2^m-1$, the other values less than $10^6$
given by these expressions are: 5, 21, 51, 85, 341, 455, 819, 1365, 3855, 5461, 13107, 21845, 29127, 31775,
87381, 209715, 258111, 349525, and 986895.

\subsubsection{Wacky Trigonometric Identities}

The computer algebra systems {\em Mathematica 5.0} and {\em Maple 7.0} will not automatically simplify the
expressions\footnote{The CASs remain sadly unreliable. For example, {\em Mathematica 5.0} computes the limit
    \[\lim_{n\to\infty} \frac{\pi n /3^n}{\sin(\pi n /3^n)} \left(1+\frac{n-1}{n}\right)\]
to be 0. The authors are unaware of any source which documents the mathematical failures of these often-used
rarely-cited closed-source programs.}
    $$\frac 1{\sin(\pi/7)}-\frac{1}{\sin(2\pi/7)}-\frac{1}{\sin (3\pi/7)}$$
    \begin{multline*}
        \frac 1{\sin(\pi/21)}-\frac 1{\sin(2\pi/21)}- \frac 1{\sin(4\pi/21)}
        - \frac 1{\sin(5\pi/21)}- \frac 1{\sin(8\pi/21)} \\
        + \frac 1{\sin(10\pi/21)}
    \end{multline*}
to 0, although they can be coaxed into verifying these identities by first replacing $\sin x$ with $\frac
{1}{2i} (e^{ix}-e^{-ix})$. But even this algebraification does not enable the CA systems to verify
    $$-2 = \sum_{k=1}^{11} \frac{\sin(2^{k+4} \pi/89)}{\sin(2^k \pi/89)}.$$
We will use the construction of perfect covers and Theorem~\ref{thm:MainTheorem} (with some further
manipulation) to generate these and other trigonometric identities.

If $q\ge 3$ and $2^m \equiv 1 \pmod{q}$, then (from the proof of Lemma~\ref{lem:constructions})
    $$0=\sum_{k=1}^m \frac{\omega^{2^{k-1}}}{1-\omega^{2^k}} = \sum_{k=0}^{m-1}
    \frac{1}{\omega^{2^k}-\omega^{-2^k}} = \sum_{k=0}^{m-1} \frac{1}{2i \sin(2^k \,2\pi/q)}.$$
Thus,
    $$0=\sum_{k=1}^m \frac{1}{\sin(2^k\pi/q)} = \sum_{k=1}^m \csc(2^k \pi /q).$$
Two of the identities above are given by this with $(q,m)=(7,3)$ and $(q,m)=(21,6)$ (using
$\sin(x)=\sin(\pi-x)=\sin(x+2\pi)$). We note that Jager \& Lenstra~\cite{MR52:3071} showed that all linear
dependencies over $\QQ$ of $\csc(\pi/q), \csc(2\pi/q), \dots, \csc(\frac{q-1}{2} \pi/q)$ have this form when $q$
is a prime.

Let us take a closer look at our basic identity, which we can
expand as follows:
    \begin{multline*}
    0=\sum_{k=0}^{m-1} \frac{1}{\omega^{2^k} - \omega^{2^{-k}}} = \\
    -\sum_{k=0}^{m-1} \left(\omega^{2^k} + \omega^{3 \cdot 2^k} + \omega^{5 \cdot 2^k} + \ldots + \omega^{(2t-1)
    \cdot 2^k} - \frac{\omega^{2t \cdot 2^k}}{\omega^{2^k} - \omega^{-2^k}}\right),
    \end{multline*}
for any $t$ with $0<t\le m-1$. Thus,
    $$\sum_{k=0}^{m-1} \frac{\omega^{2t \cdot 2^k}}{\omega^{2^k}-\omega^{2^{-k}}} =
        \sum_{k=0}^{m-1} \sum_{u=1}^t \omega^{(2u-1)2^k}.$$
Let $C_q(x)$ denote the (multi-)set $\{x \cdot 2^j \bmod{q} \colon 0  \le j <m \}$. Thus, interchanging order of
summation, we have
    \[ \sum_{k=0}^{m-1} \frac{\omega^{2t \cdot 2^k}}{\omega^{2^k}
    -\omega^{2^{-k}}} = \sum_{u=1}^t \sum_{a \in C_q(2u-1)} \omega^a.\]

We can rewrite the LHS summand as
\[  \frac{\omega^{2t \cdot 2^k}}{\omega^{2^k}
-\omega^{2^{-k}}} = \frac{\cos(2 \pi 2t \cdot\ 2^k/q) + i\sin(2 \pi 2t \cdot 2^k/q)}{2i \sin(2 \pi 2^k/q)} \]
where we have taken $\omega = e^{2 \pi i/q}.$
Thus,
    \begin{align*}
    \sum_{k=0}^{m-1} \frac{\cos(2 \pi 2t \cdot 2^k/q)}{\sin(2 \pi \cdot 2^k/q)}
        &= -2 \Im \left(\sum_{u=1}^t \sum_{a \in C_q(2u-1)} \omega ^a\right),\\
    \sum_{k=0}^{m-1} \frac{\sin(2 \pi 2t \cdot 2^k/q)}{\sin(2 \pi \cdot 2^k/q)}
        &= \phantom{-}2\Re \left(\sum_{u=1}^t \sum_{a \in C_q(2u-1)} \omega ^a\right).
    \end{align*}
Hence, we need to understand the sums
\[ S(q,t):= \sum_{u=1}^t \,\sum_{a \in C_q(2u-1)} \omega ^a.\]
To begin with, if $-1 \in C_q(1)$ then every term $\omega ^a$ in
$S$ will also have its conjugate $\omega ^{-a}$ in $S$ as well.
This happens exactly when the order of $2$ modulo $q$ is even, and
in this case we have
\[ \sum_{k=0}^{m-1} \frac{\cos(2 \pi 2t \cdot 2^k/q)}{\sin(2 \pi
\cdot 2^k/q)} = 0, \]
but for the trivial reason that every term in the sum occurs with
its negative! In general, each $C_q(2u-1)$ contains $r$ terms, so that
$S(q,t)$ is the sum of $t$ blocks of $r$ powers
 of $\omega$. If the union
of these blocks is a perfect $c$-covering of $\{\omega, \omega^2,
\omega^3, \ldots , \omega^{q-1}\}$ then $S(q,t)$ is just equal to $-c$.
For example, for $q=7$, we have $C_7(1) = \{1,2,4\},
C_7(3) = \{3,5,6\}$. Thus,
\[ S(7,2) = \sum_{u=1}^2 \sum_{a \in C_7(2u-1)} \omega^a =
(\omega + \omega^2 +\omega^4) + (\omega^3 + \omega^5 + \omega^6)
= -1, \] so that
\[ \sum_{k=0}^2 \frac{\sin(8 \pi \cdot 2^k/7)}{\sin(2 \pi \cdot 2^k/7)}
= -2 .\] Another simple case where this happens is when $q \equiv 1 \pmod{6}$ is prime, 2 has order $(q-1)/3$
modulo $q$, and 1, 3 and 5 are in distinct $C_q(i)$. In this case, $C_q(1), C_q(3)$ and $C_q(5)$ are disjoint,
so that their union is $\ZZ_q \backslash \{0\}$, which implies that $S_q(3) = -1$. Examples of this occur for
$q=229, 277, 283$, etc. Note that the real part of $\omega^a$ is equal to the real part of $\omega^{-a}$. Hence,
if the $C_q(b_i), 1 \le i \le (q-1)/m,$ form a complete set of disjoint $C's$, then in forming a perfect
covering of $\ZZ_q \backslash \{0\}$, we an use either $C_q(b_i)$ or $C_q(-b_i)$ interchangeably, if we only
want to control the real part of $S(q,t)$. For example, for $q=89$, 2 has order 11 modulo 89, and the complete
set of disjoint $C's$ is $C_{89}(1)$, $C_{89}(3)$, $C_{89}(3)$, $C_{89}(9)$, $C_{89}(11)$,
 $C_{89}(13)$, $C_{89}(19)$, and
$C_{89}(33)$. However, one can check that $-1 \in C_{89}(11)$,
 $-3 \in C_{89}(33)$, $-5 \in C_{89}(9)$ and $-13 \in C_{89}(19)$.
Thus,
    $$
    \Re\left(\sum_{u=1}^8 \sum_{a \in C_{89}(2u-1)} \omega^a\right)
    = \Re\left(\sum_{j=1}^{88} \omega^j\right) = -1.
    $$
This implies (as usual) the unlikely identity
    \[
    \sum_{k=0}^{10} \frac{\sin(32 \pi \cdot 2^k/89)}{\sin(2 \pi \cdot 2^k/89)} = -2.
    \]
Using these ideas (and other extensions thereof), many other results
of this type can be derived but we will not pursue these here.

\subsection{Bounding $q$}

We assume in this section that the $p_k$ are distinct, $0<p_k<q$, define \mbox{$g_k:=\gcd(p_k,q)$}, assume that
$\gcd(q,p_1,\dots,p_m)=\gcd(g_1,\dots,g_m)=1$, $r_k\in\ZZ$, that the $m$ sequences $\B_{p_k,r_k}$ are a perfect
$c$-fold covering, and that there is no pair $i<j$ with $p_i+p_j=q$ (this is weaker than the hypotheses of the
CFC). Let \mbox{$g:=\min\{g_1,g_2,\dots,g_m\}$}, and let $n$ be the multiplicity of $g$ in
$\{g_1,g_2,\dots,g_m\}$.

\begin{lem} \label{lem:few_gi}
If $j\not\equiv0\pmod{q}$ is a multiple of one of $g_1, g_2, \dots, g_m$, then it is a multiple of at least
three of them.
\end{lem}

\begin{proof}
Using this $j$ in the second displayed equation in Corollary~\ref{cor:CoveringCriterion}, we have
    $$0=\sum_{k} g_k \tf{g_k|j} \frac{\omega^{-jr_k}}{1-\omega^{j\bar p_k}}.$$
Clearly this sum cannot have only one nonzero term. Suppose that it has exactly two, say $g_1|j$ and $g_2|j$.
If $g_1<g_2$, then $g_1$ is a multiple of only one of $g_1,\dots,g_m$, which cannot happen (set $j=g_1$). Thus
without loss of generality we may assume that $j=g_1=g_2$. We have
    $$0=\frac{\omega^{-jr_1}}{1-\omega^{j\bar p_1}}+\frac{\omega^{-jr_2}}{1-\omega^{j\bar p_2}}.$$
Multiply by $\omega^{jr_2}$ and clear denominators to get (setting $d=r_2-r_1$)
    $$1-\omega^{j\bar p_1}+\omega^{jd}-\omega^{j(\bar p_2+d)} = 0.$$
If four complex numbers with the same modulus sum to 0, then we can the split the four into two pairs, each of
which sums to 0.

Our first case is $1=\omega^{j\bar p_1}$, $\omega^{j(\bar p_2+d)}=\omega^{jd}$, which is the same as $j\bar p_1
\equiv 0 \pmod q, j\bar p_2+jd \equiv jd \pmod q$. It follows that $j\bar p_1 \equiv j\bar p_2 \pmod{q}$, and we
multiply this congruence by $p_1p_2$ (a multiple of $j^2$) to get
    $$j^2p_2\equiv j^2 p_1 \pmod{j^2 q}.$$
Thus $p_1 \equiv p_2 \pmod{q}$, and since $0<p_k<q$, we actually have $p_1=p_2$. This contradicts our
hypothesis that the $p_k$ are distinct.

Our second case is $1=-\omega^{jd}$, $\omega^{j\bar p_1}=-\omega^{j(\bar p_2+d)}$, which forces $q$ to be even
and which is the equivalent to $jd \equiv q/2 \pmod q$, $j\bar p_1 \equiv q/2+j \bar p_2+jd \pmod q$. It follows
that $j\bar p_1 \equiv j\bar p_2\pmod q$, which we handled above.

Our third case is $1=\omega^{j(\bar p_2+d)}$, $\omega^{j\bar p_1}=\omega^{jd}$, and this is the same as $j\bar
p_2+jd \equiv 0 \pmod q$, $j\bar p_1 \equiv jd \pmod q$. Combining these gives $j\bar p_1+j\bar p_2 \equiv 0
\pmod q$. Multiply this equation by $p_1p_2$ (a multiple of $j^2$) to get
    $$0\equiv p_1p_2(j\bar p_1+j\bar p_2) \equiv j^2p_2 + j^2 p_1 \pmod{j^2 q},$$
whence $p_1+p_2 \equiv 0 \pmod{q}$. Since $0<p_k<q$, we actually have $p_1+p_2=q$. This contradicts our
hypothesis that there is no pair of `$p$'s which sum to $q$.
\end{proof}

\begin{lem} \label{lem:q.bound}
If $n=3$, then $q\le 7g$; if $n=4$, then $q\le 17g$; if $n=5$ then $q\le 33g$; if $n=6$ then $q\le730g$, and in
general $$q\le \left( \Big(\frac{n}{e-1}+1\Big)^n+1\right) g.$$
\end{lem}

\begin{proof}
Set $j=p_k$ in Corollary~\ref{cor:CoveringCriterion}, and subtract $2g_k \omega^{-p_kr_k}({1-\omega^{p_k \bar
p_k}})^{-1}$ from both sides, to get
    $$-2 g_k \frac{\omega^{-p_kr_k}}{1-\omega^{p_k \bar p_k}}
    = \sum_{i=1}^m \big(1-2\tf{k=i}\big) \,{\tf{g_i|p_k} \,g_i}\, \frac{\omega^{-p_kr_i}}{1-\omega^{p_k\bar p_i}}.$$
Taking the absolute value of each side, and using the triangle inequality and the identity
    $$1-e^{i s} = -2i e^{is/2} \sin(s/2)$$ we get
    \begin{equation*}
    \frac{2 g_k}{\sin(\pi g_k /q)} \le \sum_{i=1}^m \frac{\tf{g_i|p_k} g_i}{|\sin(\pi p_k \bar p_i /q)|}.
    \end{equation*}
for all $k\in[m]$.

Suppose that our numbering has $g=g_k$ for $k\in [n]$. We have the inequalities for $k\in[n]$
    $$\frac{2}{\sin(\pi g /q)} \leq \sum_{i=1}^n \frac{1}{|\sin( \pi p_k \bar p_i /q)|}.$$
By replacing $q$ with $q/g$, we can assume without loss of generality that $g=1$ (the bound we find for $q$
will in truth be a bound for $q/g$). We wish to show that if $q$ is large enough, then the RHS must be small
for some choice of $k$.

Let $\|x\|$ be the distance from $x$ to the nearest multiple of $q$, and let $z$ satisfy $\sum_{i=z+1}^{z+n-1}
i^{-1} <1$. Consider the directed graph with vertices $p_1,\dots,p_n$, with an edge from $p_i$ to $p_j$ if
$\|p_j \bar p_i\| \le z$. Every finite directed graph contains either a sink (a point with no out-edges) or a
cycle. If $p_{v_1},p_{v_2},\dots,p_{v_\beta}$ is a cycle, then
    \begin{multline*}
    \|p_{v_2} \bar p_{v_1}\|\cdot \|p_{v_3} \bar p_{v_2}\|  \cdots \|p_{v_1} \bar p_{v_\beta}\|
    \equiv (\pm p_{v_1} \bar p_{v_2})(\pm p_{v_2} \bar p_{v_3}) \cdots (\pm p_{v_\beta} \bar p_{v_1})\\
    \equiv \pm 1 \pmod{q}
    \end{multline*}
and
    $$1<\|p_{v_1} \bar p_{v_2}\|\cdot \|p_{v_2} \bar p_{v_3}\|  \cdots \|p_{v_\beta} \bar p_{v_1}\|
        \le z^\beta \le z^n.$$
Therefore, $z^n\ge q-1$. If $p_k$ is a sink, then all ($1\le i \le n, i\not=k$) of $\|p_k\bar p_i\|$ are
strictly greater than $z$. Since $p_1,p_2,\dots,p_n$ are distinct and there is no solution to $p_i+p_j=q$, the
$n$ values $\|p_k \bar p_i\|$ are also distinct. Thus,
    \begin{equation}\label{equ:1}
    \frac{2}{\sin(\pi/q)} \le \sum_{i=1}^n \frac{1}{|\sin( \pi p_k \bar p_i /q)|} \le
        \frac{1}{\sin(\pi/q)}+ \sum_{i=1}^{n-1} \frac{1}{\sin(\pi(z+i)/q)}.
    \end{equation}
Using the approximation $\sin x \approx x$ for small $x$, one sees that~\eqref{equ:1} bounds $q$ above. More
precisely, if the graph has a cycle, then $q\le z^n+1$, and otherwise $q$ is bounded by~\eqref{equ:1}. For
$n=3,4,5,6$, we calculate that $z=1,2,2,3$, and consequently $q$ is at most $\max\{2,7\}$, $\max\{17,10\}$,
$\max\{33,24\}$, \linebreak $\max\{730, 3\}$, respectively.

To prove the ``in general'' statement, we need only work with $n\ge 7$. Let $z=\floor{n/(e-1)+1}$. We define the
graph as above, and handle a cycle in the same way. If $p_k$ is a sink, then~\eqref{equ:1} becomes
    $$\frac{1}{\pi/q}
    \le \frac{1}{\sin(\pi/q)}
    \le  \sum_{i=1}^{n-1} \frac{1}{\sin(\pi(z+i)/q)}
    \le \frac{z+n-1}{\sin(\pi(z+n-1)/q)} \sum_{i=1}^{n-1} \frac 1{z+i},$$
where we have used the inequalities
    \[
    x\ge \sin x \ge \frac{\sin(\pi(z+n-1)/q)}{\pi(z+n-1)/q} \; x
    \]
for $0<x<\pi(z+n-1)/q$. We note that
    $$\sum_{i=1}^{n-1} \frac 1{z+i} \le \int_{z}^{z+n-1} \frac{dx}{x} = \ln\left(\frac{z+n-1}{z}\right),$$
so that our inequality can be weakened to read
    $$1 \le \frac{\pi(z+n-1)/q}{\sin(\pi(z+n-1)/q)}\ln\left(\frac{z+n-1}{z}\right).$$
This is inconsistent for any $z>n/(e-1)$, $q>z^n$ and $n\ge 4$.
\end{proof}

\begin{lem} \label{lem:g=1} If $m\le 5$, then $g=1$.
\end{lem}

\begin{proof}
Renumber so that $g=g_1 \le g_2\le \dots \le g_m$. If all the $g_i$ are equal, then $g=1$ since
$\gcd\{g_1,\dots,g_m\}=1$. Thus we may assume that $n<m$.

By Lemma~\ref{lem:few_gi} with $j=g$, $n\ge 3$. Since $n<m$, we may assume that $4\le m\le6$, so that
$g=g_1=g_2=g_3<g_4$. Clearly $g$ divides $g_1,g_2,g_3$ by definition, and $g_4,g_5$ by Lemma~\ref{lem:few_gi}
(with $j=g_4$ and $j=g_5$, respectively). Since $\gcd\{g_1,\dots,g_m\}=1$, we have shown that $g=1$.
\end{proof}

\begin{thm}
The CFC is true for $m\le 5$.
\end{thm}

\begin{proof}
First, observe that by Lemma~\ref{lem:g=1} we may restrict our attention to sequences $\B_{p_k,r_k}$ with $p_1$
relatively prime to $q$. Moreover, as a consequence of Corollary~\ref{cor:CoveringCriterion}, the sequences
$\B_{p_1,r_1},\B_{p_2,r_2},\dots,\B_{p_m,r_m}$ are a perfect covering if and only the sequences
    $$
        \B_{1,0},\quad \B_{\bar p_1 p_2, p_1 r_2}, \quad\dots,\quad \B_{\bar p_1 p_m, p_1r_m}
    $$
are a perfect covering.\footnote{This will typically affect the {\em multiplicity} of the covering.} Thus, we
may assume that $p_1=1$ and $r_1=0$. Corollary~\ref{cor:CoveringCriterion} also tells us that $p_m\equiv -
\sum_{k=1}^{m-1} p_k \pmod{q}$, and with $j=1$ that
    \begin{equation}
    0=\sum_{k=1}^m \tf{g_k=1} \frac{\omega^{-r_k}}{1-\omega^{\bar p_k}}.
    \end{equation}
This implies that
    \begin{equation}
    \label{equ:computationinequality}
    1 \le \sum_{k=2}^m \left|\frac{1-\omega}{1-\omega^{\bar p_k}}\right| = \sum_{k=2}^m
    \frac{\sin(\pi/q)}{\sin(\pi \bar p_k/q)}.
    \end{equation}
In this expression there are $m-1$ degrees of freedom: $q$, $p_2$, $p_3$, $\dots$, $p_{m-1}$. Without loss of
generality $1=p_1<p_2<p_3<\dots<p_m<q$. Also, the hypotheses of the CFC imply that there is no solution to
$p_i+p_j=q$, and Lemma~\ref{lem:q.bound} (with Lemma~\ref{lem:g=1}) implies $q\le33$.

There are only 346 $m$-tuples $(p_1,\dots,p_m)$ with $m\le 5$, $$1=p_1<p_2<\dots<p_m<q\le33,$$ $\sum_{k=1}^m p_k
\equiv 0 \pmod{q}$, no solution to $p_i+p_j=q$, and satisfying~\eqref{equ:computationinequality}. We refine our
search by noting that if $(p_k,q)=1$, then $(\bar p_k p_1, \bar p_k p_2, \dots,\bar p_k p_m)$ must also be on
our list of 346 (this is equivalent to taking values of $j$ other than 1 in deriving the
inequality~\eqref{equ:computationinequality}). This pares the list down to a single tuple for $m=3$, a single
tuple for $m=4$, and 10 tuples for $m=5$. The tuples predicted by the CFC are on these lists, and the remaining
9 tuples are eliminated by an exhaustive search for $r_2,\dots,r_5$ such that
    $0=\frac{1}{1-\omega}+\sum_{k=2}^5 \frac{\omega^{-r_k}}{1-\omega^{\bar p_k}}.$

The only perfect coverings with 5 or fewer Beatty sequences are those predicted by the Covering Fraenkel
Conjecture.
\end{proof}

We remark that we may similarly reduce the $m=6$ case of Fraenkel's Conjecture (but not the CFC) to a finite
computation. For example, if $m=6,n=4$, then we may argue as in the $m=5$ case of Lemma~\ref{lem:g=1} that
$g=1$, and so by Lemma~\ref{lem:q.bound} we get the bound $q\le 730$. If $m=6,n=3$ (so that $q\le 7g$), then we
may (using Lemma~\ref{lem:few_gi}) renumber so that $g=g_1=g_2=g_3$ and $h=g_4=g_5=g_6$, with $\gcd(g,h)=1$.
Since the $p_k$ are distinct, we have $p_1+p_2+p_3\ge g+2g+3g=6g$ and $p_4+p_5+p_6\ge 6h>6g$. Thus $q=\sum p_k
>12 g$, a contradiction. In contrast, for $m=7$, we arrive at the consistent inequalities $q\le 17g$ and $q=\sum
p_k
>16g$.

\section{Proving Fraenkel's Conjecture}

We envision a non-computational proof of Fraenkel's Conjecture along the following lines. Suppose that
$\B_{p_k,r_k}^q$ ($1\le k \le m$) partition $\ZZ$, with $p_1<p_2<\dots<p_m<q$, and suppose that $m$ is minimal.
Let $g_k:=\gcd(p_k,q)$.

Now, suppose that \mbox{$g:=\min_k\{g_k\}$} is larger than 1. The Covering Criterion with $j=g$ yields
    $$
    0=\sum_{\substack{k=1\\g_k=g}}^m \frac{\omega^{-gv_k}}{1-\omega^{g \bar p_k}}
    $$
Let $u_1,\dots,u_n$ be those $\bar p_k$ for which $g_k=g$ (there are $n\ge3$ of them by Lemma~\ref{lem:few_gi},
and $n<m$ since $\gcd\{g_1,\dots,g_m\}=1$), and let $v_1,\dots,v_n$ be the negatives of the $r_k$ for which
$g_k=g$. Also replace $\omega=e^{2\pi i /q}$ with $x=e^{2\pi i g/q}$, and we get
    $$
    0=\sum_{k=1}^n \frac{x^{v_k}}{1-x^{u_k}}
    $$
This seems to imply that the sum vanishes for $x$ {\em any} $q/g$-th root of unity, which would imply that these
sequences alone form a perfect covering, whence $\sum_{k=1}^m p_k \tf{g_k=g} \ge q$. This is impossible since
$n<m$. Thus, the following conjecture implies that, in the present setting, $g=1$.

\begin{cnj}\label{cnj:rationalfunction}
Suppose that $1\le u_1 < u_2 < \dots < u_n <q$, with $\gcd(u_k,q)=1$ for all $k$, and let $v_1, \dots, v_n$ be
arbitrary integers. If the function
    $$f(x) := \sum_{k=1}^n \frac{x^{v_k}}{1-x^{u_k}}$$
vanishes at $x=e^{2\pi i /q}$, then $\sum_{k=1}^n u_k \ge q$.
\end{cnj}

This in turn simplifies the Covering Criterion substantially. Considering the absolute value of the Covering
Criterion as in the proof of Lemma~\ref{lem:q.bound}\footnote{In Lemma~\ref{lem:q.bound} we found that for
$q>n^n$ (roughly), $k$ could be chosen to make a particular inequality invalid.
Conjecture~\ref{cnj:StrongMartin} posits that the $n^n$ bound can be improved to $2^n$, and barring the single
exception of $q=2^n-1$, it can be improved to $(7/4)^n$. This is supported by computational investigations.} the
following conjecture becomes relevant.
\begin{cnj}\label{cnj:StrongMartin}
Suppose that $p_1,p_2,\dots, p_n$ are distinct and relatively prime to $q>(7/4)^n$, with $\sum p_k \le q$, and
for each $k\in[n]$
   $$\frac{2}{\sin(\pi /q)} \leq \sum_{i=1}^n \frac{1}{|\sin( \pi p_k \bar p_i /q)|}.$$
Then $q=2^n-1$ and $\{p_1,\dots, p_n\}\equiv \{1,2,\dots, 2^{n-1}\} \pmod{q}$.
\end{cnj}

At this point, we would have shown that a counterexample to Fraenkel's Conjecture (with $m$ sequences) must have
$q<(7/4)^m$. We envision handling this situation combinatorially, probably in conjunction with Tijdeman's
combinatorial restrictions. He~\cite{MR2001f:11039}*{Lemma 4} notes that if $\B_{p_k,r_k}$ and $\B_{p_j,r_j}$
are disjoint, then either $p_k=p_j$ or $\floor{q/p_k}\not=\floor{q/p_j}$, and so
    $$\floor{q/p_m} < \dots < \floor{q/p_2} < \floor{q/p_1}.$$
Also, his main lemma~\cite{MR2001f:11039}*{Lemma 3} can be strengthened (using the same proof, but in terms of
Beatty sequences instead of balanced sequences) to provide the powerful restriction on the $p_k$'s in a
counterexample with minimal $m$: $p_k \le (q-2g_k)/3.$ Tijdeman used these two lemmas (and some casework) to
show that $m\ge7$. Thus, the remaining situation has many sequences with small (but spread out) $p_k$'s and
quite small $q$.

\section{Refining the Conjectures}

Ideally, one would like arithmetic conditions on $\alpha_k,r_k$ for the sequences $\{\floor{n\alpha_k+r_k}\colon
n\in \ZZ\}$ to be a perfect covering, without assuming that the $\alpha$ are distinct. Morikawa has given such
conditions for a small number of sequences to be a perfect 1-cover, see~\cite{MR2001h:11011} for a brief
description of Morikawa's work and citations for his many papers on the topic.

The Covering Fraenkel Conjecture that we have advanced is another step in this direction. A more ambitious step
would be to replace the condition ``with no proper subset $I\subsetneq[m]$ having $\sum_{i\in I} p_i\equiv 0
\pmod{q}$'' condition with the weaker condition ``with no proper subset of the sequences being a perfect
covering.''

We note that Conjecture~\ref{cnj:rationalfunction} can likely be strengthened:
\begin{cnj}
Suppose that $1\le u_1 < u_2 < \dots < u_n <q$, with $\gcd(u_k,q)=1$ for all $k$, and with no subset of the
$u_k$'s summing to a multiple of $q$, and let $v_1, \dots, v_n$ be arbitrary integers. If the function
    $$f(x) := \sum_{k=1}^n \frac{x^{v_k}}{1-x^{u_k}}$$
vanishes at $x=e^{2\pi i /q}$, then it vanishes at all $q$-th roots of unity except $x=1$.
\end{cnj}

Joe Buhler notes that
    $$\frac{1}{1-x}+\frac{x^5}{1-x^2}+\frac{x^{10}}{1-x^4}+\frac{x^{10}}{1-x^{11}}
        +\frac{x^5}{1-x^{13}}+\frac{1}{1-x^{14}}$$
vanishes at primitive 15-th roots of unity, but not at the primitive 5-th roots of unity, and thus the condition
on sums of subsets of the $u_k$'s is necessary.

\section*{Acknowledgements}
The authors wish to thank the National Science Foundation for supporting mathematics in general and the authors
in particular. We thank Professor Greg Martin of the University of British Columbia for supplying the main idea
underlying Lemma~\ref{lem:q.bound}.

\end{document}